# DE LA K-THÉORIE ALGÉBRIQUE VERS LA K-THÉORIE HERMITIENNE

Max KAROUBI
Mathématiques, UMR 7586 du CNRS, case 7012, Université Paris 7, 2, place Jussieu, 75251 Paris cedex 05, France
e.mail : karoubi@math.jussieu.fr

**FROM ALGEBRAIC K-THEORY TO HERMITIAN K-THEORY**

In this paper, we introduce a new morphism between Algebraic and Hermitian K-theory. The topological analog is the Adams operation $\psi^2$ in real K-theory. From this morphism, we deduce a lower bound for the higher Algebraic K-theory of a ring A in terms of the classical Witt group of the ring $A \otimes A^{op}$.

## 1. La construction fondamentale.

Soit A une k-algèbre unitaire, telle que $1/2 \in A$. Nous nous proposons de construire une "classe caractéristique" nouvelle

$$K(A) \longrightarrow L(A \otimes A^{op})$$

Ici $K(A)$ désigne la K-théorie algébrique de A et $L(A \otimes A^{op})$ la K-théorie hermitienne de l'anneau $B = A \otimes_k A^{op}$, muni de l'antiinvolution permutant les deux facteurs.

Plus généralement, cette construction algébrique permet de décrire dans le § 3 une application continue entre les espaces classifiants des K-théories algébrique et hermitienne

$$\mathbf{K}(A) \longrightarrow \mathbf{L}(B)$$

d'où on déduit une application

$$K_n(A) \longrightarrow L_n(B)$$

qui est un homomorphisme de groupes si $n > 0$.

La construction est basée sur le foncteur suivant. Soit E un A-module à droite projectif de type fini et soit $E^* = \mathrm{Hom}_A(E, A)$ son dual, qui est un A-module à gauche par la règle

$$(\lambda . f)(x) = \lambda(f(x))$$

où $f \in E^*$ et où les lettres grecques désignent des éléments de A. Alors $E^* \otimes_k E$ est un B-module à droite par la règle du produit suivante

$$(f \otimes y)(\lambda \otimes \mu) = \mu f \otimes y \lambda$$

Le B-module F peut être muni de la forme hermitienne $\varphi$ suivante

$$\varphi [(f \otimes y), f' \otimes y')] = f(y') \otimes f'(y)$$

Si, pour simplifier, on pose $u = f \otimes y$ et $u' = f' \otimes y'$ et $b = \lambda \otimes \mu$, on a en effet les identités suivantes

$$\varphi(ub, u') = \varphi[(\mu f \otimes y\lambda), f' \otimes y')] = \mu f(y') \otimes f'(y)\lambda = (\mu \otimes \lambda)(f(y') \otimes f'(y))$$
$$= \overline{b}\, \varphi(u, u')$$
$$\varphi(u, u'b) = \varphi(f \otimes y, \mu f' \otimes y'\lambda) = f(y')\lambda \otimes \mu f'(y) = (f(y') \otimes f'(y))(\lambda \otimes \mu)$$
$$= \varphi(u, u')b$$

et enfin $\varphi(u', u) = f'(y) \otimes f(y') = \overline{\varphi(u, u')}$.

Par ailleurs, la forme hermitienne $\varphi$ permet d'identifier F à son antidual F* (en tant que B-module à droite), grâce au morphisme antilinéaire $\theta : F \longrightarrow F^*$ défini par $\theta(u)(u') = \varphi(u, u')$. En effet, par additivité, il suffit de le vérifier pour $E = A$, ce qui est évident.

## 2. Interprétation en termes de traces.

Supposons que A soit un anneau commutatif avec involution, notée $a \mapsto \overline{a}$. On définit un homomorphisme $A \otimes A^{op} \longrightarrow A$ par la formule $a \otimes b \mapsto a\overline{b}$. Par extension des scalaires, le produit tensoriel $(E^* \otimes E) \otimes_{A \otimes A^{op}} A$ s'identifie à $\mathrm{Hom}_A(E, \overline{E})$ [où $\overline{E}$ désigne le A-module conjugué de E]. Plus précisément, cette identification $\phi$ est définie par

$$(f \otimes y) \otimes a \mapsto [x \mapsto \overline{f(x)}\, y\, a]$$

On vérifie aisément que $\phi([\mu f \otimes y \lambda] \otimes a) = \phi([f \otimes y] \otimes \lambda \overline{\mu} a)$

**PROPOSITION**. *Compte tenu des identifications précédentes, la forme hermitienne sur* $\mathrm{Hom}_A(E, \overline{E})$ *est définie par la formule suivante* :
$$(f, g) \mapsto \mathrm{Tr}(\overline{f}\, g)$$
*où $\overline{f}$ est l'application linéaire de $\overline{E}$ dans E sous-jacente à f.*

**Exemples** : si A est le corps des complexes **C** et $E = \mathbf{C}$, $\mathrm{Hom}_A(E, \overline{E})$ s'identifie à **C** et la forme hermitienne $\phi$ n'est autre que la forme standard $(\lambda, \mu) \mapsto \overline{\lambda}\mu$. D'une manière générale, la proposition précédente permet d'associer à tout A-module E une **famille** de formes hermitiennes $\phi_\sigma$, $\sigma$ parcourant l'ensemble des involutions de l'anneau A.

## 3. Retour à la K-théorie.

Soit E est un A-module projectif de type fini et soit de nouveau $B = A \otimes A^{op}$. On note $\psi(E)$ le B-module projectif de type fini $E^* \otimes E$, muni de la forme hermitienne définie dans le § 1,

**PROPOSITION.** *Soient* E *et* F *deux* A-*modules projectifs de type fini. On a alors un isomorphisme canonique de* B-*modules hermitiens*

$$\psi(E \oplus F) = \psi(E) \oplus \psi(F) \oplus H(E^* \otimes F)$$

*où* H *est le foncteur hyperbolique standard.*

*Démonstration*. Notons d'abord que $E^* \otimes F$ est bien un B-module par la règle de multiplication suivante (déjà explicitée plus haut si E = F) :
$$(u \otimes v)(\lambda \otimes \mu) = \mu u \otimes v \lambda$$
Par ailleurs, on a un isomorphisme
$$(E \oplus F)^* \otimes (E \oplus F) \cong (E^* \otimes E) \oplus (F^* \otimes F) \oplus [(E^* \otimes F) \oplus (F^* \otimes E)]$$
Enfin, les B-modules $E^* \otimes F$ et $F^* \otimes E$ sont en dualité par l'accouplement induit
$$[(f \otimes y), (f' \otimes y')] \mapsto f(y') \otimes f'(y)$$
analogue à celui défini plus haut (si E = F). Il s'en suit que $[(E^* \otimes F) \oplus (F^* \otimes E)]$, muni de la forme induite, est isomorphe au module hyperbolique $H(E^* \otimes F)$.

Pour étendre $\psi$ au groupe de Grothendieck K(A), on remarque que la correspondance $\psi$ vérifie l'identité suivante sur le semi-groupe des classes d'isomorphie de B-modules hermitiens
$$\psi(x + y) = \psi(x) + \psi(y) + \gamma(x, y)$$

où $\gamma$ est une fonction biadditive (c'est le foncteur hyperbolique appliqué au produit tensoriel du dual de x par y). D'après une remarque due à Dold [3], on peut étendre la fonction $\psi$ au groupe de Grothendieck en posant
$$\psi(x - y) = \psi(x) - \psi(y) - \gamma(x, y) + \gamma(y, y)$$
On vérifie que cette formule a un sens : $\psi(x - y)$ ne dépend bien que de la classe de x - y dans la K-théorie de A. En d'autres termes, elle définit une application (pas un homomorphisme)
$$K(A) \longrightarrow L(B)$$

En particulier, si y est la classe d'un module libre $A^n$, on voit que

$$\psi(E - A^n) = \psi(E) - \psi(A^n) + H((A^n)^* \otimes A^n) - H(E^* \otimes A^n)$$

Il convient de noter que le foncteur $E \mapsto \psi(E)$ donne naissance à une représentation remarquable (pour $E = A^n$)
$$GL_n(A) \longrightarrow O_{p,q}(A \otimes A^{op})$$
avec $p = n + (n^2 - n)/2$ et $q = (n^2 - n)/2$. Cependant, il n'est pas évident à priori d'en déduire une application continue
$$BGL(A)^+ \longrightarrow BO(A \otimes A^{op})^+$$
entre les espaces classifiants usuels des K-théories algébrique et hermitienne

Pour définir une telle application (et expliciter ses propriétés multiplicatives), nous allons procéder de manière indirecte en interprétant la K-théorie algébrique (et hermitienne) en termes de fibrés plats. Par exemple, si X est un CW-complexe fini, la K-théorie algébrique de X à coefficients dans A peut être définie en termes de fibrés plats "virtuels" sur X, dont la fibre est un A-module projectif de type fini (cf. [5] pour les détails). Plus précisément, le groupe de Grothendieck des classes d'isomorphie de tels fibrés, noté $K_A(X)$, s'identifie à l'ensemble des classes d'homotopie de X dans l'espace $K(A) \times BGL(A)^+$ où K(A) est muni de la topologie discrète. En particulier, les groupes $K_n(A)$ peuvent être définis en termes de fibrés plats sur des sphères homologiques de dimension n.

On définit de manière analogue la K-théorie hermitienne de X à coefficients dans un anneau R avec antinvolution, notée $L_R(X)$, qui s'identifie aussi à l'ensemble des classes d'homotopie de X dans $L(R) \times BO(R)^+$. Puisque la correspondance $E \mapsto \psi(E)$ est fonctorielle, elle définit une **application** entre les groupes de Grothendieck de fibrés virtuels associés, par la méthode de [3] avec R = B, soit

$$K_A(X) \longrightarrow L_B(X)$$

On en déduit une application continue bien définie à homotopie faible près entre les espaces classifiants correspondants. En choisissant pour X une sphère de dimension n > 0 par exemple, on a un homomorphisme de groupes

$$K_n(A) \longrightarrow L_n(A \otimes A^{op})$$

qui étend l'application définie précédemment pour n = 0.

## 4. Interprétation topologique

Supposons maintenant que A = C(X) soit l'algèbre des fonctions continues à valeurs réelles sur un espace compact X (on pourrait aussi considérer les fonctions à valeurs complexes, en munissant A de l'involution définie par la conjugaison). Le complété naturel de $A \otimes A^{op} = A \otimes A$ est alors l'algèbre des fonctions continues sur X x X, munie de l'involution associée à la permutation des variables. Le choix d'une métrique sur les fibrés permet d'identifier la K-théorie hermitienne de C(X x X) à la K-théorie équivariante $K_{\mathbf{Z}/2}(X \times X)$. L'analogue topologique de l'homomorphisme précédent $\psi$ est alors la "power operation" d'Atiyah [1]

$$K(X) \longrightarrow K_{\mathbf{Z}/2}(X \times X)$$

où K(X) désigne la K-théorie topologique réelle. En particulier, si on se restreint à la diagonale de X x X, on trouve une opération $\psi$ réduite

$$\psi' : K(X) \longrightarrow K_{\mathbf{Z}/2}(X) \cong K(X) \oplus K(X)$$

dont les deux composantes sont les opérations de Grothendieck $S^2(x)$ et $\lambda^2(x)$. La différence $S^2(x) - \lambda^2(x)$ est l'opération d'Adams $\psi^2$ considérée comme morphisme d'anneaux de K(A) dans l'anneau de Witt W(A) qui s'identifie à K(A) d'après un théorème bien connu. On notera que $\psi^2$ est la multiplication par une puissance de 2 sur la K-théorie réduite des sphères. Par conséquent, c'est un isomorphisme modulo la 2-torsion si X est un CW-complexe fini. En effet, le foncteur de K-théorie **topologique** transforme une suite exacte courte d'algèbres de Banach en une longue suite exacte de groupes de K-théorie topologique.

## 5. Un exemple d'application à la K-théorie algébrique supérieure.

Puisque l'opération $\psi$ s'étend aux espaces classifiants, elle induit un homomorphisme sur les théories à coefficients (avec $B = A \otimes A^{op}$)

$$K_n(A ; \mathbf{Z}/p^\alpha) \longrightarrow W_n(B ; \mathbf{Z}/p^\alpha)$$

pour $n \geq 0$. Ceci est bien défini si $n = 0$ en convenant de poser $K_0(A ; \mathbf{Z}/p^\alpha) = K_0(A)/p^\alpha$ et $W_0(B ; \mathbf{Z}/p^\alpha) = W_0(B)/p^\alpha$.

Supposons maintenant que p soit un nombre premier impair et désignons par $H_n$ l'image de $K_n(A ; \mathbf{Z}/p^\alpha)$ dans $W_n(B ; \mathbf{Z}/p^\alpha)$. Nous comptons trouver une minoration de $H_n$, donc de $K_n(A; \mathbf{Z}/p^\alpha)$, en utilisant le théorème de périodicité en K-théorie hermitienne [4].

Soit $\mathbf{Z}' = \mathbf{Z}[1/2]$. D'après [2] il existe un "élément de Bott" $b_{4m}$ dans $K_{4m}(\mathbf{Z}'; \mathbf{Z}/p^\alpha)$ dont l'image dans le groupe de K-théorie topologique $K_{4m}^{top}(\mathbf{R}, \mathbf{Z}/p^\alpha) \cong \mathbf{Z}/p^\alpha$ est le générateur topologique (avec $m = p^{\alpha-1}(p-1)/2$).

**THEOREME**. *Avec les notations précédentes, on a un diagramme commutatif* [1]

$$\begin{array}{ccc} K_n(A ; \mathbf{Z}/p^\alpha) & \longrightarrow & W_n(B ; \mathbf{Z}/p^\alpha) \\ \downarrow \gamma & & \downarrow \beta \\ K_{n+4m}(A ; \mathbf{Z}/p^\alpha) & \longrightarrow & W_{n+4m}(B ; \mathbf{Z}/p^\alpha) \end{array}$$

*où $\gamma$ est le cup-produit par $b_{4m}$ et où $\beta$ est l'homomorphisme de périodicité en K-théorie hermitienne. Puisque $\beta$ est un isomorphisme pour $n > 0$ et injectif pour $n = 0$, on en déduit que l'ordre du groupe $H_{n+4m}$ est minoré par celui de $H_n$.*

*Démonstration*. L'image du générateur de Bott $b_{4m}$ dans $W_{4n}(\mathbf{R} ; \mathbf{Z}/p^\alpha)$ est un générateur de ce groupe qui s'identifie lui-même à $W_{4m}^{top}(\mathbf{R} ; \mathbf{Z}/p^\alpha) \cong \mathbf{Z}/p^\alpha$ d'après [4]. On conclut la démonstration en remarquant que les morphismes horizontaux sont compatibles avec les structures multiplicatives en K-théorie algébrique et hermitienne.

---

[1] On note $W_n$ le groupe de Witt "supérieur" qui est le conoyau de la flèche naturelle de $K_n$ dans $L_n$, induite par le foncteur hyperbolique.

# Références